%% file: paper.tex
\documentclass{article}
\input{myhead}

\begin{document}

\title{Fourier analysis of the CGMN method for solving the Helmholtz equation}
\author{Tristan van Leeuwen,\\ {\small Centrum Wiskunde \& Informatica, Amsterdam, The Netherlands}}
\maketitle
\begin{abstract}
The Helmholtz equation arises in many applications, such as seismic and medical imaging.
These application are characterized by the need to propagate many
wavelengths through an inhomogeneous medium. The typical size of the problems in 3D applications
precludes the use of direct factorization to solve the equation and hence iterative methods are used in practice.
For higher wavenumbers, the system becomes increasingly indefinite and thus good preconditioners
need to be constructed. 
In this note we consider an accelerated Kazcmarz method (CGMN) 
and present an expression for the resulting iteration matrix.
This iteration matrix can be used to analyze the convergence of the CGMN method.
In particular, we present a Fourier analysis for the method applied to the 1D Helmholtz equation.
This analysis suggests an optimal choice of the relaxation parameter.
Finally, we present some numerical experiments.
\end{abstract}

\section{Introduction}
The Helmholtz equation arises in many applications such as seismic and medical imaging.
These application are characterized by the need to propagate many
wavelengths through an inhomogeneous medium. The scale of the problems in 3D application precludes the use
of direct factorization to solve the equation and hence iterative methods are used in practice.
For higher wavenumbers, the system becomes increasingly indefinite and thus good preconditioners
need to be constructed. For an overview of preconditioning techniques for the Helmholtz equation 
we refer to \cite{ernst11} and references cited therein.

In this note, we consider an accelerated Kazcmarz method (CGMN, \cite{Bjorck1979}) 
for solving the Helmholtz equation \cite{Gordon2013}. We present an expression of 
the iteration matrix that allows us to analyze its convergence behaviour.

We proceeds as follows.
First, we give a brief overview of the CGMN method and its relation to SSOR. 
This gives us an expression for the iteration matrix in terms of the original matrix.
Then, we present the Fourier
analysis for a 1D finite-difference discretization of the Helmholtz equation and 
we present some numerical examples. Finally, we draw conclusions and discuss future work.

\subsection{The CGMN method}

The CGMN method \cite{Bjorck1979} relies on a symmetric Kaczmarz (SKACZ) sweep in conjunction with the conjugate 
gradient (CG) method. 

The Kaczmarz method solves a system of $N$ equations, $A\mbf{u} = \mbf{s}$, by
cyclically projecting the iterate onto rows of the matrix \cite{Kaczmarz1937}
\bq
\mbf{u} := \mbf{u} + \omega\left(s_i - \transp{\mbf{a}_i}\mathbf{u}\right)\mbf{a}_i/||\mbf{a}_i||_2^2,\quad i = 1\ldots N,
\eq
where $\mbf{a}_i$ denotes the $i$-th row of $A$ as column vector and $0<\omega< 2$ is a relaxation parameter. Note that 
we may also let $\omega$ vary per row.
Introducing the matrices $Q_i = \left(I - \omega\mbf{a}_i\transp{\mbf{a}_i}/||\mbf{a}_i||_2^2 \right)$,
we may write this iteration as
\bq
\mbf{u}:=  Q_i\mathbf{u} + \omega\mbf{a}_i/||\mbf{a}_i||_2^2
\eq
A \emph{double sweep} through the matrix (from row 1 to $N$ and back) can then be denoted by
\bq
\mbf{u}:= Q\mbf{u} + R\mbf{s},
\eq
where $Q = Q_0Q_1\ldots Q_{N-1}Q_{N-1}\ldots Q_0$ and $R$ contains all the
factors multiplying $\mbf{b}$.
It is easily verified that the matrices $Q_i$ are symmetric and have rank 1 with eigenvalue $1-\omega$.
It follows that $Q$ is symmetric and has eigenvalues $\in [-1,1]$. Hence, $(I-Q)$ is symmetric and positive semi-definite
and we can use CG to solve the equivalent system 
\bq
(I - Q)\mbf{u} = R\mbf{s}.
\eq
The question is, how does this new system behave in terms of $\omega$? In the next section, we will
first derive a convenient expression for the matrix $Q$ in terms of the matrix $A$, and in the 
subsequent section we present a Fourier analysis of the error propagation applied to the 1D Helmholtz equation.

\subsection{Relation to SSOR}

It is well-known that a SKACZ sweep on the system $A\mbf{u}=\mbf{s}$ is equivalent
to a SSOR sweep on the normal equations $AA^T\mbf{x} = \mbf{s}$, $\mbf{u}=A^T\mbf{x}$ \cite{Bjorck1979}. 
Splitting the matrix $AA^T$ into its diagonal, strictly lower and upper triangular parts $D, L$ and $L^T$
such that $AA^T = D + L + L^T$,  the SSOR iteration can be represented as follows \cite{saad}
\bq
\mbf{x}:=G\mbf{x} + H\mbf{s},
\eq
where
\bq
\label{eq:G}
G &=& (D+\omega L^T)^{-1}((1-\omega)D - \omega L)(D+\omega L)^{-1}((1-\omega)D - \omega L^T),\\
\label{eq:H}
H &=& \omega(2-\omega)(D+\omega L^T)^{-1}D(D+\omega L)^{-1}.
\eq
Since the SKACZ and SSOR iterations are equivalent, the matrices $Q,R$ and $G,H$ are related via
\bq
QA^T &=& A^TG,\\
R    &=& A^TH.
\eq
Next, we note that we may write $G$ in terms of $H$ and $A$ as
\bq
G &=& I - HAA^T,
\eq
from which it follows that we can write
\bq
\label{eq:Q}
Q &=& I - A^THA.
\eq

\section{Fourier analysis of the error propagation}
We study the properties of the preconditioner by looking at the \emph{error propagation matrix}.
This matrix relates errors in $\mbf{u}$ when doing a standard Richardson iteration on 
the preconditioned system:
\bq
\mbf{u}_{k+1} = (I - (I-Q))\mbf{u}_{k+1} + R\mbf{b}.
\eq
In terms of the error $\mbf{e}_{k+1} = \mbf{u}_{k} - \mbf{u}^*$, we get
\bq
\mbf{e}_{k+1} = Q\mbf{e}_{k},
\eq
Similarly, the residual $\mbf{r}_{k} = A\mbf{e}_k$ is propagated by
\bq
\mbf{r}_{k+1} = G^T\mbf{r}_k.
\eq

We analyze the behaviour of the SKACZ preconditioner through a Fourier analysis of the error
propagation matrix \cite{brandt77, kettler}. We decompose the error into its Fourier modes $e_j(\theta)=e^{\imath j \theta}$
and want to find the amplitude $a(\theta)$ such that
\bq
Q\mbf{e}(\theta) = a(\theta)\mbf{e}(\theta), 
\eq
Using eq. (\ref{eq:H}) and (\ref{eq:Q}) we can factorize the amplitude function in terms of those of the matrices 
$A$ ($a_1$), $(D+\omega L)$ ($a_2$), $D$ ($a_3$) and $(D+\omega L^T)$ ($a_4$):
\bq
\label{eq:amp}
a = 1 - \omega(2-\omega) \frac{a_1^2 a_3}{a_2a_4}.
\eq
Note that the matrix $G$ has the same amplitude function.

\subsection{1D case}
We consider a simple finite difference discretization of the the 1D Helmholtz equation $(k^2 + \mathrm{d}^2/\mathrm{d}x^2)u = s$
on the domain $\Omega = [0,1]$ with Dirichlet boundary conditions. The resulting system of equations is denoted by 
\bq
A\mbf{u} = \mbf{s},
\eq
where 
\bq
a_{ii} &=& k^2 - 2h^{-2},\\
a_{i,i\pm1} &=& h^{-2},
\eq
and $h=1/(N+1)$ denotes the gridspacing and $N$ denotes the number of interior gridpoints. 
The matrix $A$ is symmetric and has eigenvectors
\bq
\mbf{v}_i^n = \sin(n\pi i h), \quad  n = 1,2, \ldots, N
\eq
with corresponding eigenvalues
\bq
\lambda_n = k^2 + 2h^{-2}(\cos(n\pi h) - 1).
\eq

The entries of the matrices $D$ and $L$ are given by
$d_{ii} = ||\mbf{a}_i||_2^2 =  \gamma^2 + 2h^{-4}$, 
$l_{i,i-1} = \mbf{a}_i^T\mbf{a}_{i-1} = 2h^{-2}\gamma$ and $l_{i,i-2}=h^{-4}$ where $\gamma = k^2-2h^{-2}$.

It is now readily verified that
\bq
a_1 &=& \gamma + 2h^{-2}\cos\theta,\\
a_2 &=& \gamma^2 + 2h^{-4} + \omega h^{-2}\bigl(2\gamma\exp^{-\imath\theta} + h^{-2}\exp^{-2\imath\theta} \bigr),\\
a_3 &=& \gamma^2 + 2h^{-4},\\
a_4 &=& \gamma^2 + 2h^{-4 } + \omega h^{-2}\bigl(2\gamma\exp^{+\imath\theta} + h^{-2}\exp^{+2\imath\theta} \bigr).
\eq

Substituting these expressions in eq. (\ref{eq:amp}) yields 
\bq
\label{eq:a1d}
\lefteqn{a(\theta,\omega,h) = 1 - } \nonumber\\
&&\frac{\omega(2-\omega)\beta\bigl(\gamma+2h^{-2}\cos\theta\bigr)^2}
{\beta^2+2\beta\omega(2\gamma h^{-2}\cos\theta+h^{-4}\cos 2\theta)+\omega^2h^{-4}(4\gamma^2+4\gamma h^{-2}\cos\theta+h^{-4})},\nonumber\\
\eq
where $\beta =  \gamma^2+2h^{-4}$.

\begin{figure}
\centering
\begin{tabular}{ccc}
\includegraphics[scale=.2]{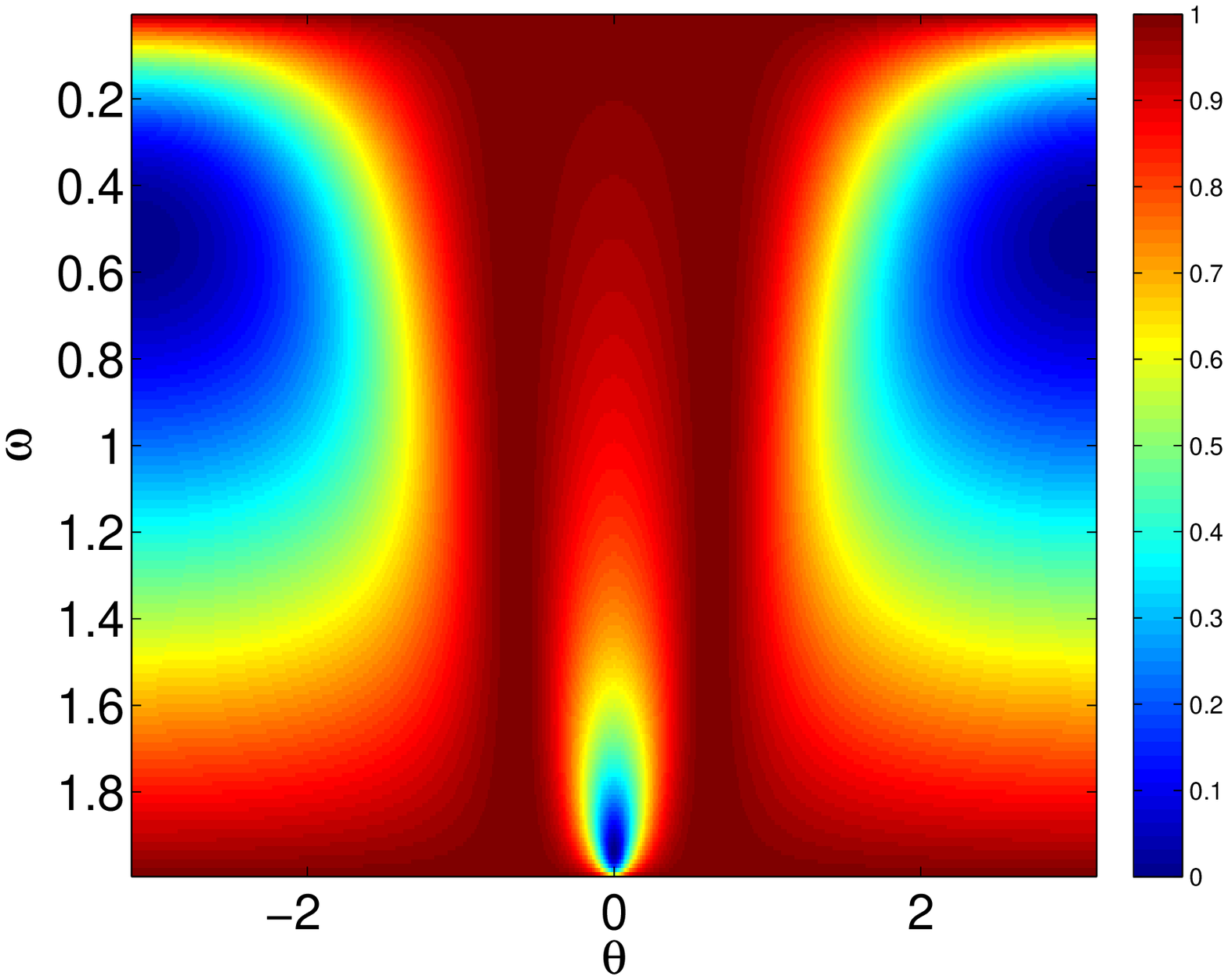}&
\includegraphics[scale=.2]{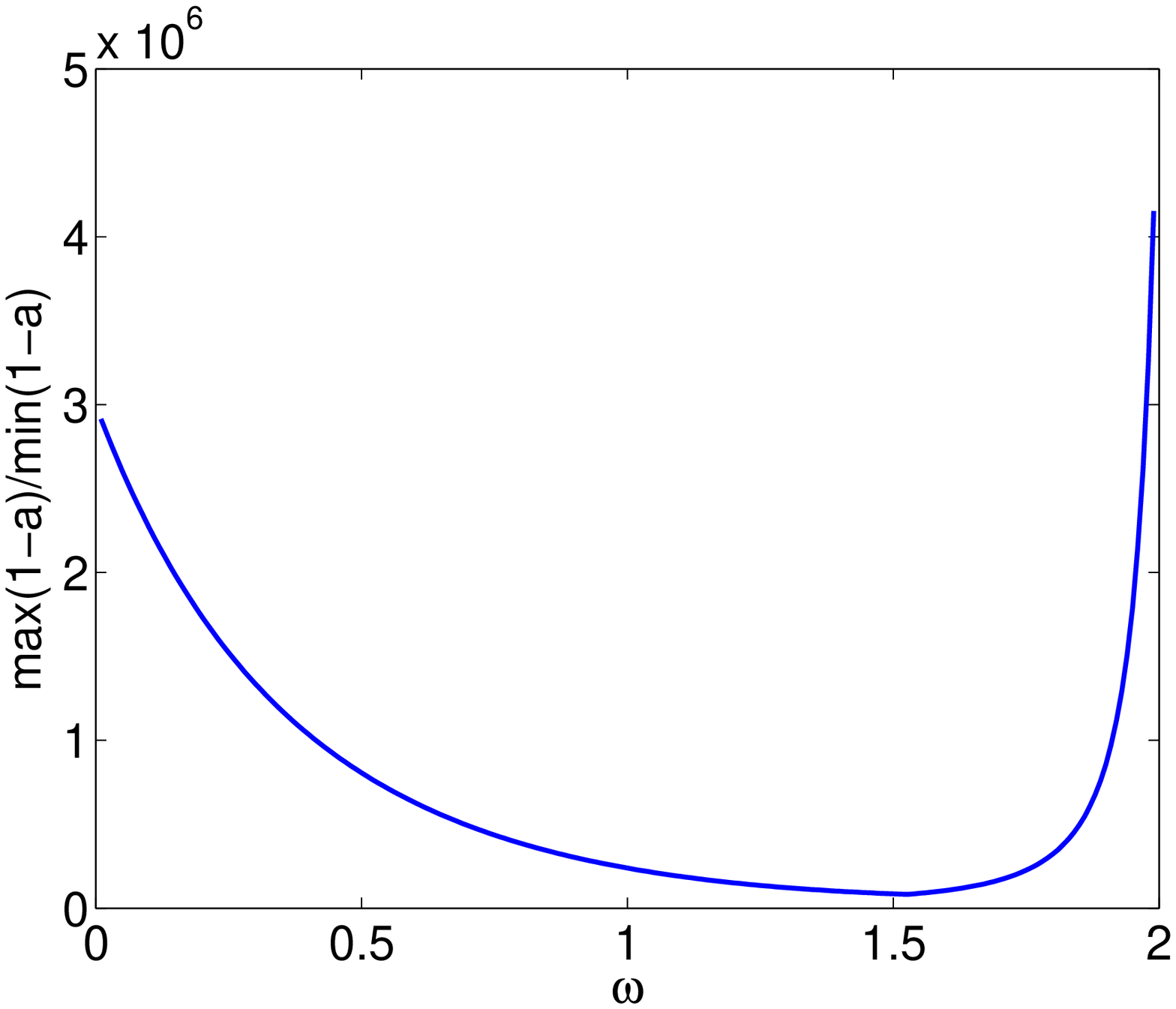}&
\includegraphics[scale=.2]{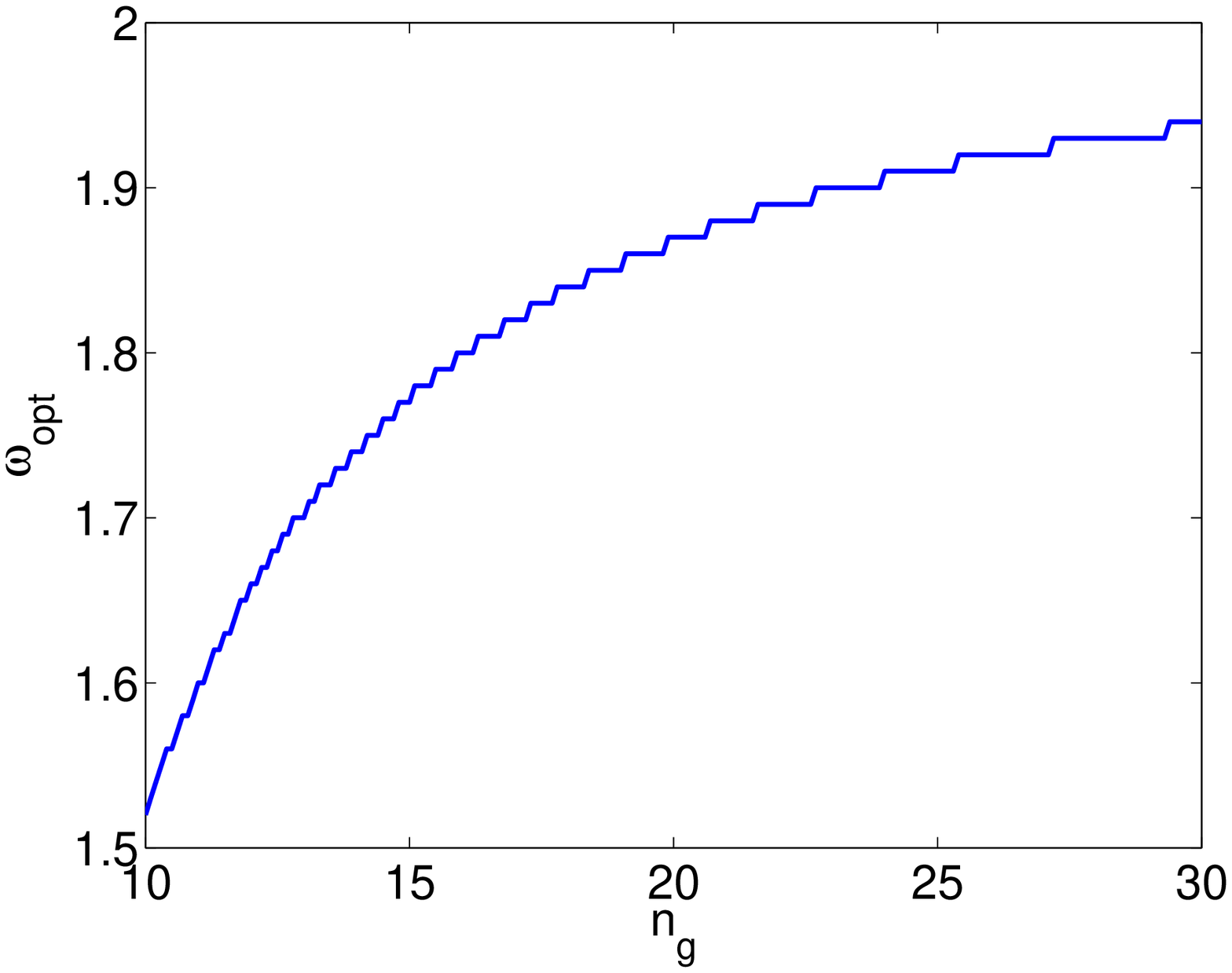}\\
{\small (a)}&{\small (b)}&{\small (c)}\\
\end{tabular}
\label{fig:error}
\caption{(a) Amplitude (cf. eq. \ref{eq:a1d}) as a function of $\theta$ and $\omega$ for a fixed number of gridpoints
per wavelength $n_g=10$. (b) Shows the `condition number' as a function of $\omega$ for $n_g=10$. The optimal $\omega \approx 1.5$.  (c) shows the optimal $\omega$ as a function of $n_g$.}
\end{figure}

To ensure a certain amount of gridpoints per wavelength, $n_g$, we let $kh = 2\pi/n_g$. Note that in 
this case, the amplitude only depends on $n_g$ and $\omega$. The question is, what is the optimal $\omega$
for a given $n_g$. Figure \ref{fig:error} (a) shows the amplitude as a function of $\theta$ and $\omega$ for $n_g=10$.
In Figure \ref{fig:error} (b) we show $\frac{\max_{\theta}(1-a(\theta))}{\min_{\theta}(1-a(\theta))}$, which is roughly equivalent to 
plotting the condition number of $I-Q$, and this suggests an optimal value of $\omega\approx 1.5$ for $n_g=10$.
Figure \ref{fig:error} (c) shows the optimal $\omega$ as a function of $n_g$. 
We can use this curve to adapt $\omega$ to the (local) wavenumber.

\section{Numerical experiments}
We present some numerical experiments that verify the results discussed above.

\subsection{1D}
\begin{figure}
\centering
\begin{tabular}{cc}
\includegraphics[scale=.3]{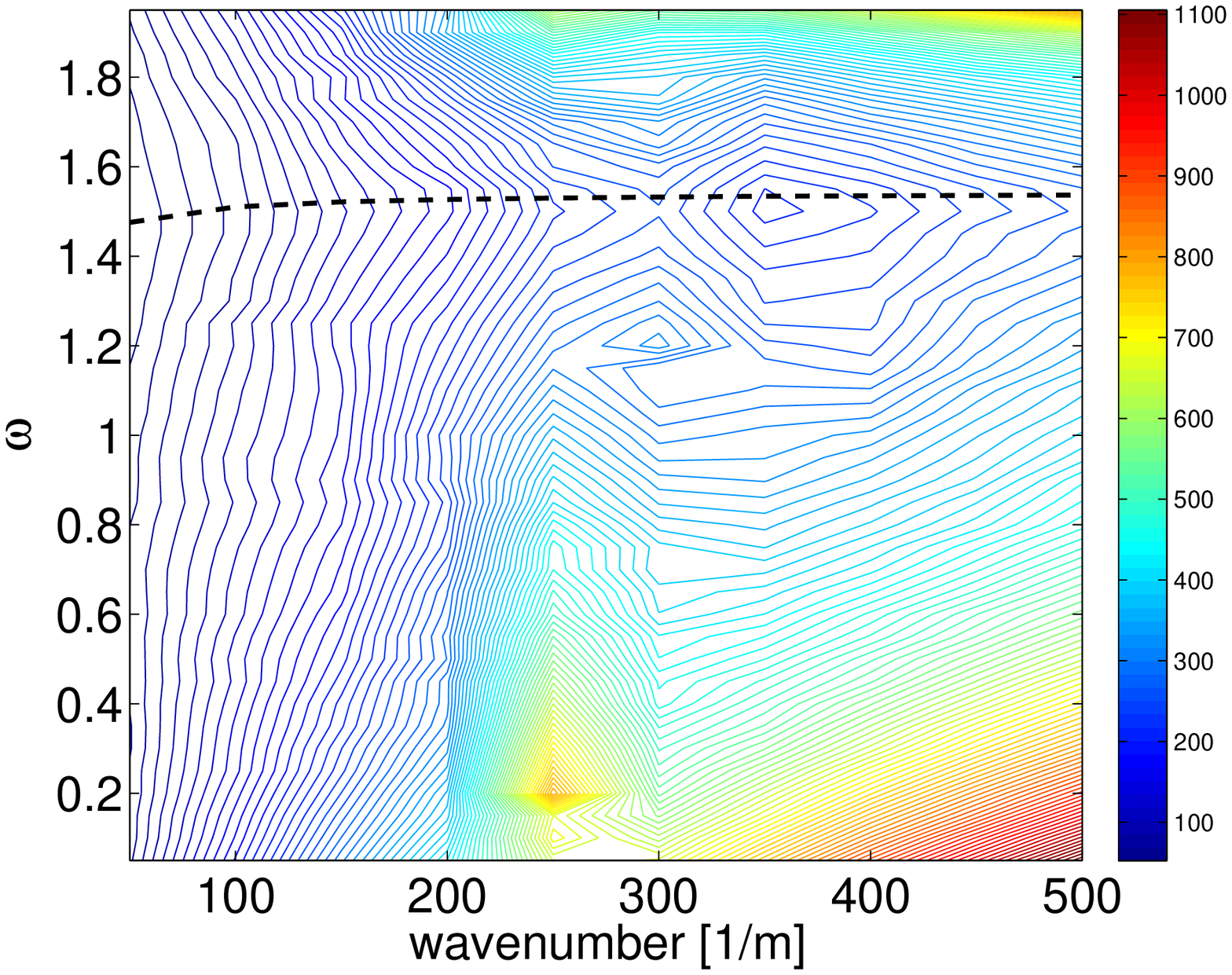}&
\includegraphics[scale=.3]{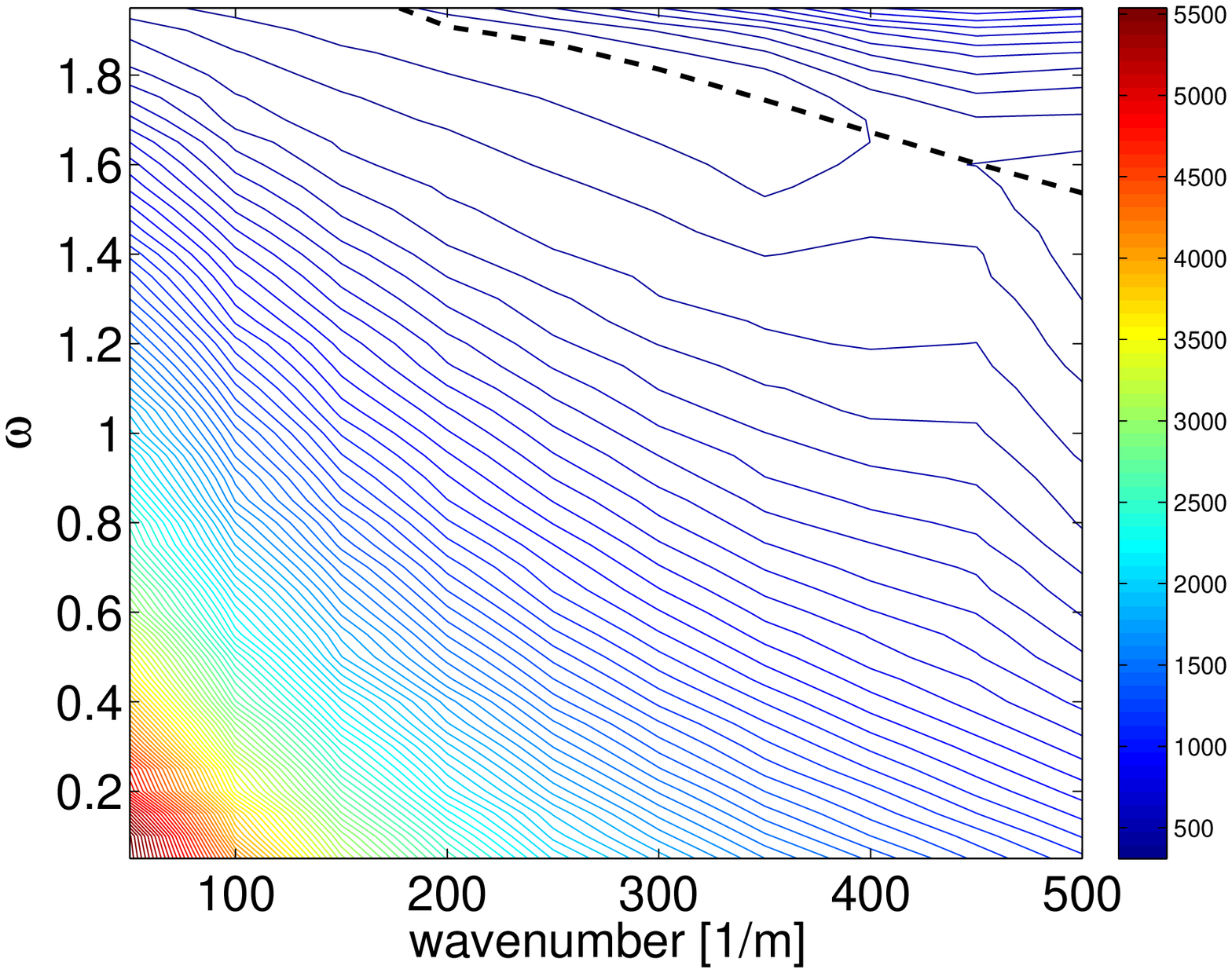}\\
{\small (a)}&{\small (b)}\\
\end{tabular}
\label{fig:exp1d}
\caption{Number of CG iterations needed to converge to a tolerance of $10^{-6}$ for various 
$\omega$ with either a fixed number of gridpoints per wavelength $n_g=10$ (a) or a fixed gridspacing based on the highest wavenumber used (b).
The predicted optimal $\omega$ is indicated by a dashed line and coincides with the lowest iteration count.}
\end{figure}

In the first experiment
we solve the Helmholtz equation for various wavenumbers $k$, with a fixed number of $n_g=10$ gridpoints per wavelength. 
Figure \ref{fig:exp1d} (a) shows the number of CG iterations needed to converge to a tolerance of $10^{-6}$ for various 
$\omega$. The predicted optimal value (cf. figure \ref{fig:error} (c)) is also indicated and nicely coincides with 
the optimal $\omega$ in terms of the iteration count.

In the second experiment, we use the same values of $k$, but choose the gridspacing
based on the \emph{highest} wavenumber used. This way, the number of gridpoints per
wavelength is much higher for the lower wavenumbers. The result is shown in Figure \ref{fig:exp1d} (b). 
Again, the predicted optimal value of $\omega$ is close to the empirical optimum.

\subsection{2D}
\begin{figure}
\centering
\begin{tabular}{cc}
\includegraphics[scale=.3]{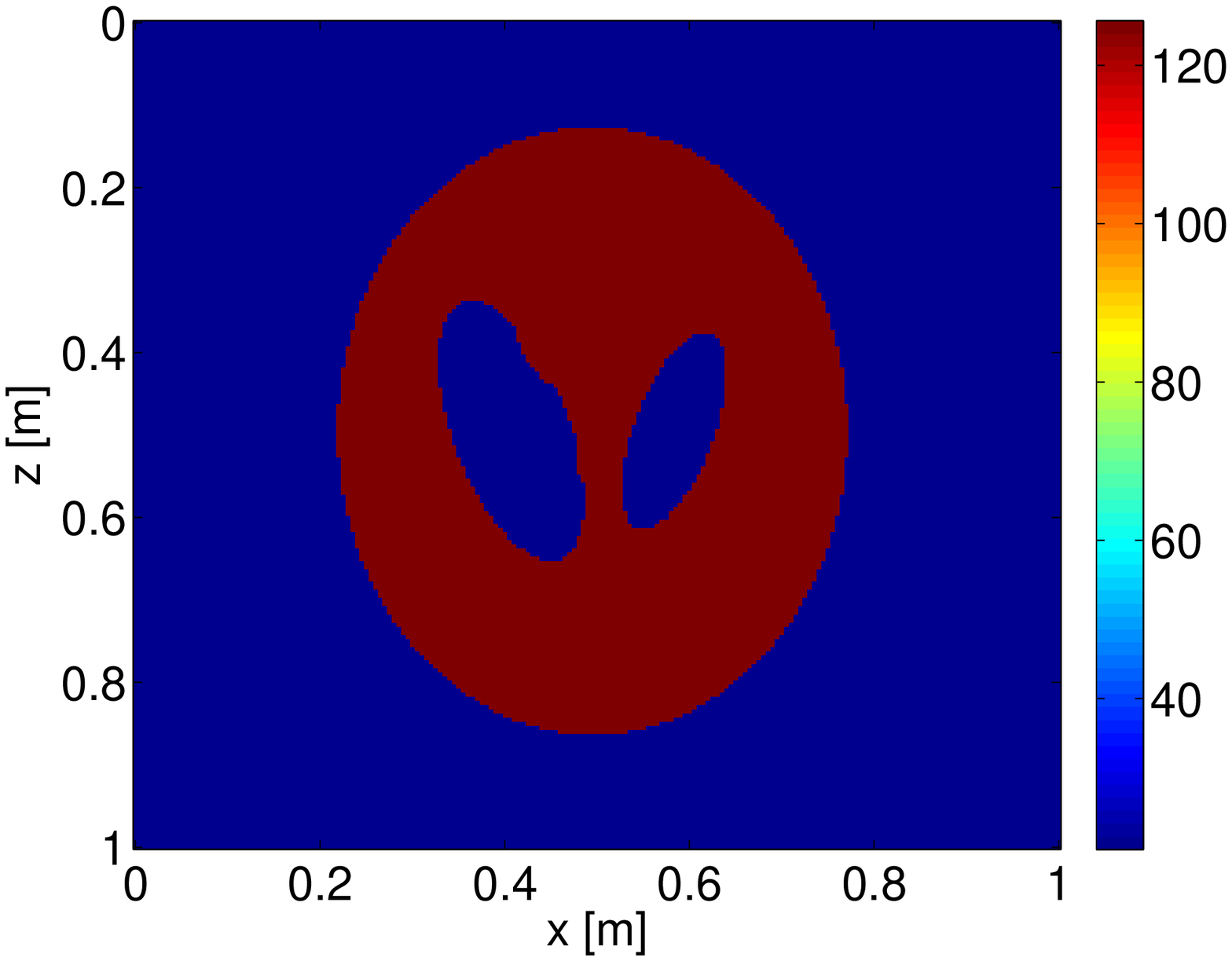}&
\includegraphics[scale=.3]{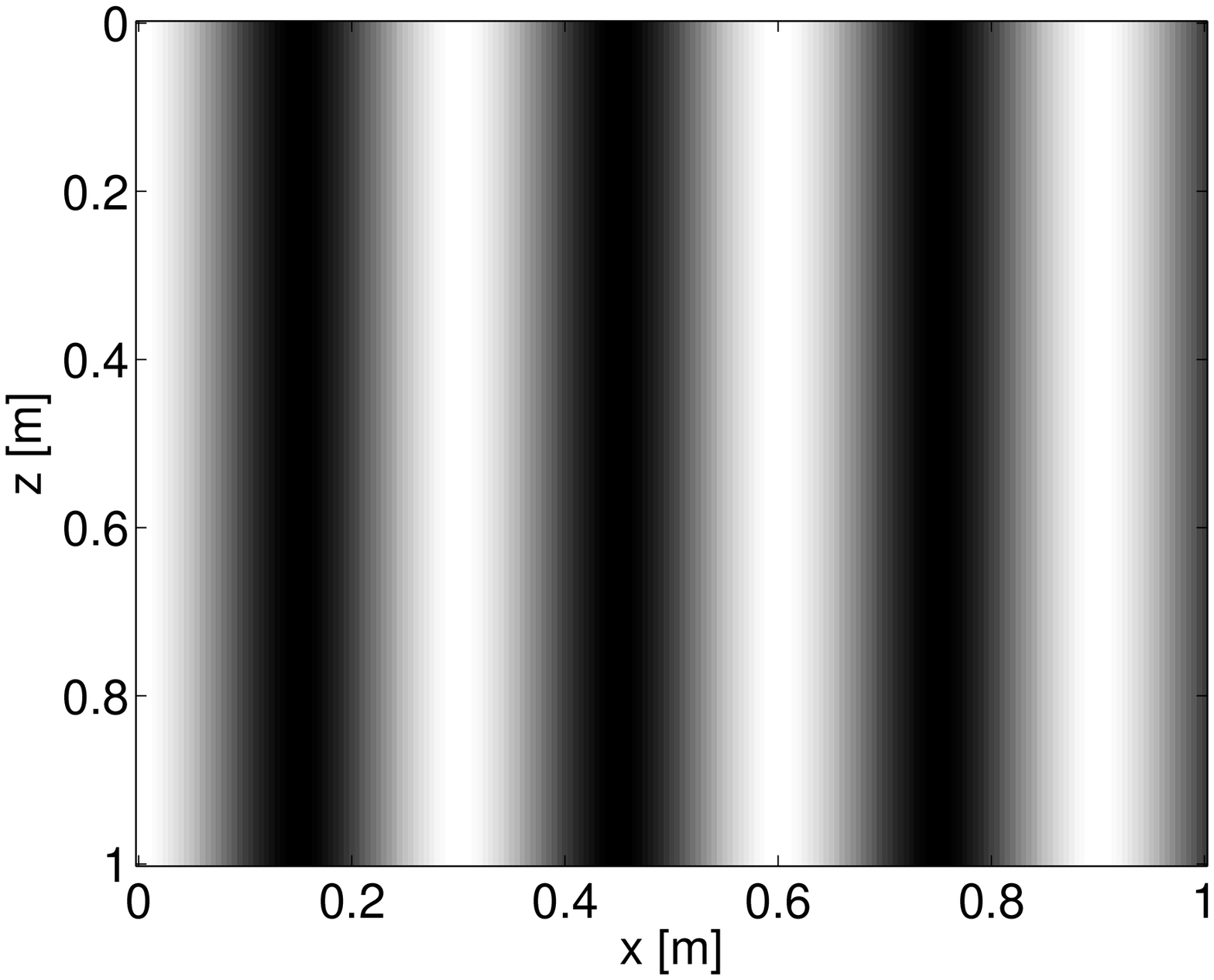}\\
{\small (a)}&{\small (b)}\\
\includegraphics[scale=.3]{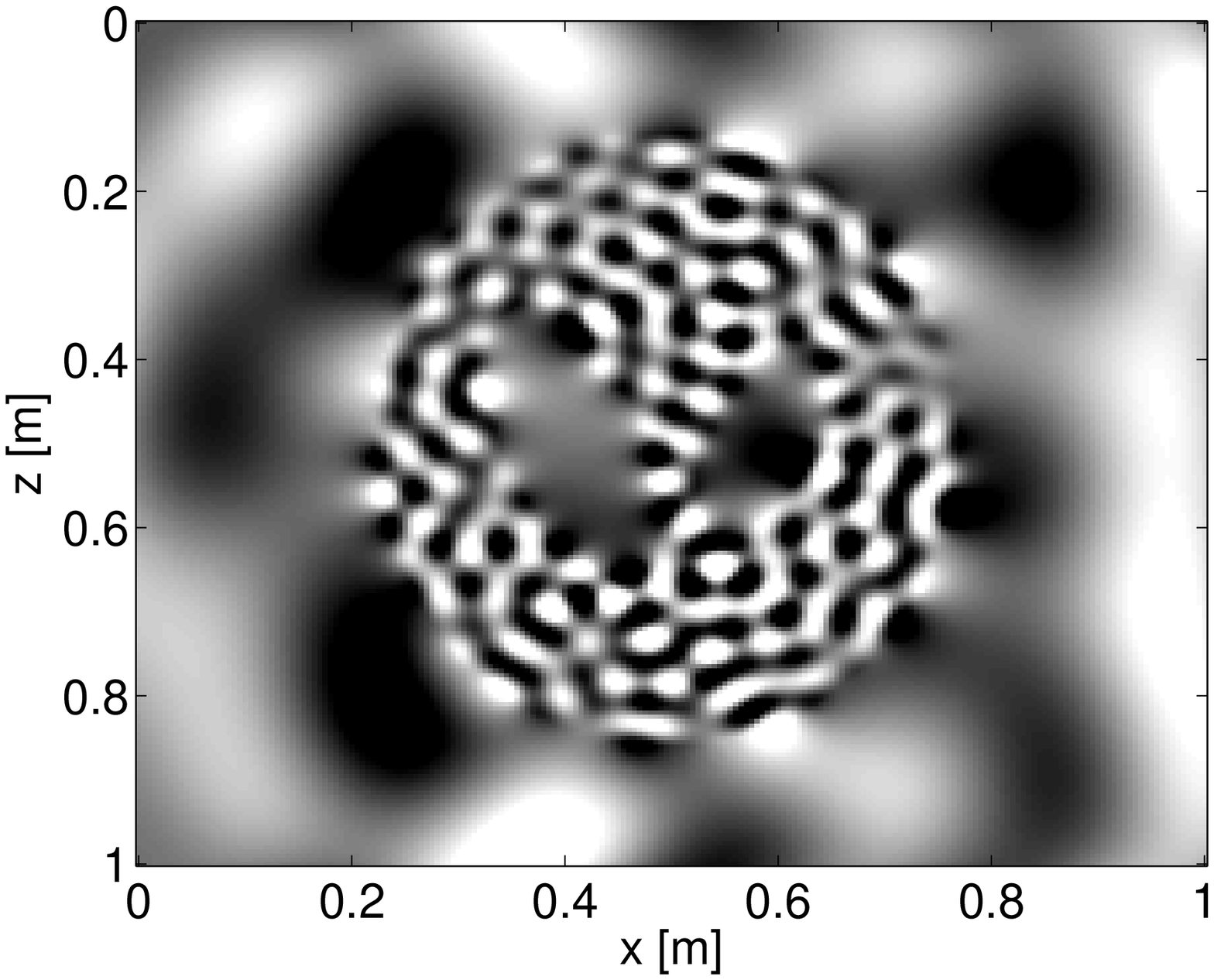}&
\includegraphics[scale=.3]{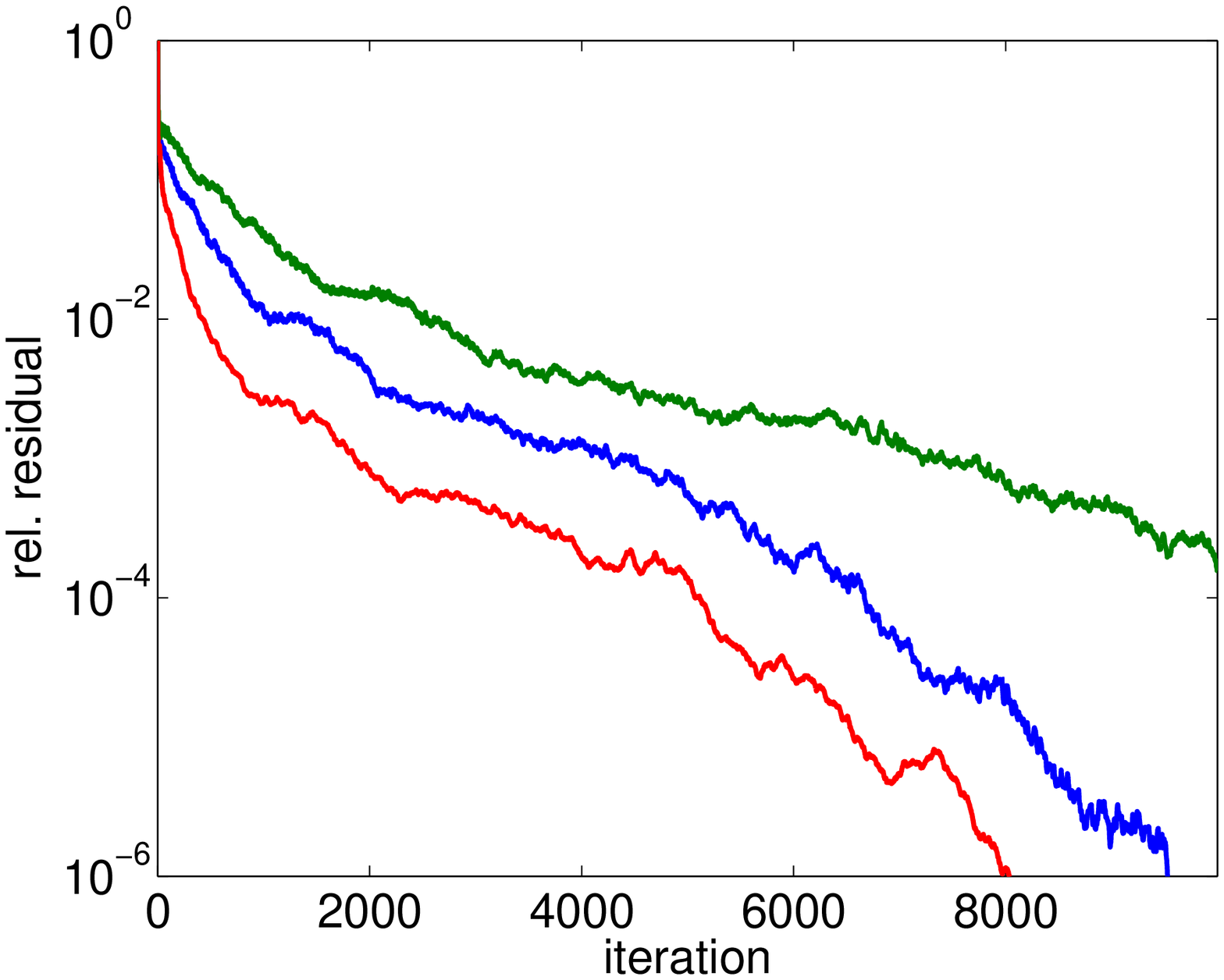}\\
{\small (c)}&{\small (d)}\\
\end{tabular}
\label{fig:exp2d}
\caption{(a) wavenumber profile, (b) incident plane wave, (c) resulting wavefield and (d) convergence histories for
$\omega = 1.5$ (blue), $\omega=1.95$ (green) and a spatially varying $\omega$ (red).}
\end{figure}

We use a medium that consist of a high-wavenumber anomaly embedded in a background with constant
wavenumber $k_0$, as depicted in figure \ref{fig:exp2d} (a). The wavefield resulting
from an incident planewave $\exp(\imath k_0 x)$ (figure \ref{fig:exp2d} (b)) is depicted in figure \ref{fig:exp2d} (c). 
The convergence histories using either a constant $\omega=1.5$ (chosen according to the highest wavenumber)
 or $\omega=1.95$ (chosen according to the lowest wavenumber) or a spatially varying $\omega$ (corresponding to the local wavenumber)
 are shown in figure \ref{fig:exp2d} (c). 
It appears that for a medium with a large contrast it may indeed be beneficial to vary the relaxation parameter
spatially.

\section{Conclusion and future work}
We have presented a Fourier analysis of the CGMN method applied to the 1D Helmholtz equation.
The amplitude function of the error propagation matrix
suggests that the optimal relation paramater depends on the number of gridpoints per wavelength 
$n_g=2\pi/(kh)$. This observation was confirmed by counting the number of CGMN iterations needed to actually solve
the system for various values of $\omega$ and $n_g$. The same optimal relation can be used to 
adapt the relaxation parameter to the local wavenumber in case of inhomogeneous media. A numerical example
on a 2D medium with a very high contrast suggests that this might indeed improve convergence.

Future research is aimed at extending this analysis to
2D/3D Helmholtz equations and varying media, as well as the parallel extension of the CGMN method presented by \cite{Gordon2013}.
An optimal strategy to adapt the relaxation parameter to the medium and perhaps
vary it per iteration is of particular interest.

\clearpage
\section*{Acknowledgments}
The author thanks Dan and Rachel Gordon for valuable discussions on this topic.
This work was in part financially supported by the Natural Sciences and 
Engineering Research Council of Canada Discovery Grant (22R81254) and the 
Collaborative Research and Development Grant DNOISE II (375142-08). 
This research was carried with support from the sponsors of the SINBAD consortium.

\bibliographystyle{plain}
\bibliography{mybib}

\clearpage

\end{document}

%% file: myhead.tex
\usepackage{amsmath}
\usepackage{amssymb}
\usepackage{amsfonts}
\usepackage{epsfig}
\usepackage{graphicx}
\usepackage{multicol}
\usepackage{color}
\usepackage{crop}
\usepackage{float}
\usepackage{booktabs}

\newcommand{\mbf}[1]{\mathbf{#1}}

\newcommand{\transp}[1]{{#1}^{\small T}}

\newcommand{\bq}{\begin{eqnarray}}
\newcommand{\eq}{\end{eqnarray}}

\renewcommand{\eqref}[1]{equation (\ref{eq:#1})}

\newcommand{\qed}{\nobreak \ifvmode \relax \else
      \ifdim\lastskip<1.5em \hskip-\lastskip
      \hskip1.5em plus0em minus0.5em \fi \nobreak
      \vrule height0.75em width0.5em depth0.25em\fi}
\renewcommand{\Pi}{\mbox{\LARGE$\pi$}}